\documentclass[a4paper,12pt]{amsart} \usepackage{euscript,eufrak,verbatim}
\usepackage[final]{epsfig}
\usepackage{subfigure}

\def\Z{\mathbb{Z}} \def\R{\mathbb{R}} \def\S{\mathbb{S}} 

\def\T{\mathbb{T}}
\def\rp{R$^+$-poor }
\def\op{O$^+$-poor }
\newenvironment{proofof}[1]{\noindent {\bf Proof of #1.}}{ \hfill\qed\\ }
\renewenvironment{proof}{\noindent {\bf Proof.}}{ \hfill\qed\\ }

\def\k{\vec{k}}
\def\l{\vec{l}}
\def\bk{\bar{k}}
\def\bl{\bar{l}}
\def\hk{\hat{k}}

\def\bg{\bar{g}}
\def\hg{\hat{g}}
\def\g{\vec{g}}

\def\phi{\varphi}
\def\be{\begin{equation}}
\def\ee{\end{equation}}

\def\bea{\begin{eqnarray}}
\def\eea{\end{eqnarray}}

\newtheorem{theorem}{Theorem}[section]
\newtheorem{lemma}[theorem]{Lemma}
\newtheorem{sublemma}[theorem]{Sublemma}
\newtheorem{corollary}[theorem]{Corollary}
\newtheorem{condition}[theorem]{Condition}
\newtheorem{convention}[theorem]{Convention}
\newtheorem{remark}[theorem]{Remark}
\newtheorem{definition}[theorem]{Definition}
\DeclareMathSymbol{\varnothing}{\mathord}{AMSb}{"3F} \begin{document}

\title{Ergodicity of Two Hard Balls in Integrable Polygons}

\author{P\'eter B\'alint} \author{Serge Troubetzkoy}

\address{Centre de Physique Th\'eorique and\\
Federation de Recherches des Unites de Mathematique de Marseille\\
CNRS Luminy, Case 907, F-13288 
Marseille Cedex 9, France\\ 
on leave from Alfr\'ed R\'enyi Institute of Mathematics 
of the H.A.S., H-1053, Re\'altanoda u. 13--15, Budapest, Hungary\\ 
and from Mathematical Institute, Technical University of Budapest, 
H-1111 Egry J\'ozsef u. 1., Budapest, Hungary}
\email{balint@cpt.univ-mrs.fr, bp@renyi.hu}
\urladdr{http://www.renyi.hu/{\lower.7ex\hbox{\~{}}}bp/}

\address{Centre de Physique Th\'eorique\\
Federation de Recherches des U\-ni\-tes de Mathematique de Marseille\\
Institut de math\'ematiques de Luminy and\\
Universit\'e de la M\'editerran\'ee\\ 
Luminy, Case 907, F-13288 Marseille Cedex 9, France}

\email{troubetz@iml.univ-mrs.fr}
\urladdr{http://iml.univ-mrs.fr/{\lower.7ex\hbox{\~{}}}troubetz/} \date{}
\begin{abstract}
We prove the hyperbolicity, ergodicity and thus the Bernoulli property
of two hard balls in one of the following four polygons:
the square, the equilateral triangle, the $45-45-90^\circ$ triangle 
or
the $30-60-90^\circ$ triangle.
\end{abstract}
\maketitle

\pagestyle{myheadings}

\markboth{ERGODICITY OF TWO BALLS}{P\'ETER B\'ALINT AND SERGE
TROUBETZKOY}

\section{Introduction}\label{sec1}
\setcounter{equation}{0}

In the past thirty years there has been much research on the 
Sinai-Boltzmann
ergodic hypothesis:  {\em the motion of an arbitrary number of elastic 
hard balls on the $\nu$-dimensional torus ($\nu \ge 2$) is ergodic} 
\cite{Doma-survey}. 
Recent results by Sim\'anyi prove this hypothesis for almost all geometric
parameters (i.e.~masses) \cite{nan1,nan2}.
However there is only one article on physically more realistic containers, 
in \cite {nandbox} Sim\'anyi has
proven the ergodicity of two 
balls in the $\nu$-dimensional cube ($\nu \ge 2$).
In this article we consider the two dimensional case and  extend this 
result 
to several other containers.
Let $P$ be a square, an equilateral triangle, a $45-45-90^\circ$ triangle 
or
a $30-60-90^\circ$ triangle.  Such polygons are called {\it integrable} since
the billiard motion of an elastic point particle reduces to an integrable system, 
the linear flow on a flat torus \cite{MT}.  
We consider the billiard system of two hard balls with unit mass and
equal radius $r$ moving uniformly in $P$ 
and colliding elastically with each other and at the
boundary $\partial P$ of $P$. Our main result is that for $r$
sufficiently small (such that the phase space is connected) 
this billiard system is hyperbolic 
(Corollary~\ref{noname}) and ergodic (Theorem~\ref{main}). By standard
argumentation, the system is then also K-mixing and Bernoulli 
(Remark~\ref{noname2}). We have specialized to the two dimensional case
for clarity of exposition, however our method allows to prove the
ergodicity of the two ball system in $n+2$-dimensional right prisms of the form 
$P \times [0,1]^n$,
where $P$ is one of the above polygons (see Subsection~\ref{noname3}).
A trivial application of our
method also allows to prove the hyperbolicity and ergodicity of one or
two balls of a specific radius
in some other tables, including $C^\infty$ convex examples and even fractal
ones (see Subsection~\ref{noname4}).

The basic idea of the proof is to lift the system to a cylindric billiard
on $\T^4$.  This is the part of the proof which is sensitive to the choice
of the polygon $P$, if we start with another rational polygon (i.e.\ all
angles are rational multiples of $180^\circ$), then we can construct this lift, but
the manifold has a more complicated structure (see Subsection \ref{noname5}).
It is the product of two copies of a higher genus flat surface with 
singularities. 

The group $G$ of symmetries generated by reflections in the sides of the  polygon
$P$ plays an important role.  For the square this group has the following special
properties which do not hold for the other polygons considered: 1) $G$ is commutative,
2) each $g \in G$ is an idempote and 3) any two invariant sets are orthogonal.
Note that for groups generated by reflections these three properties are equivalent.
The structure of $G$ is reflected in the fact that in the case of the square the
upstairs dynamical system is an orthogonal 
cylindric billiard in the sense of \cite{S}, cf. 
Remark~\ref{square}.
Furthermore, unlike Sim\'anyi, we do not use the ergodicity 
of the two ball system on  $\T^2$ \cite{nandbox}.  
We are not aware of any reasonable notion of center
of mass for higher genus surfaces \cite{Gal}, 
thus we feel that our proof is a better starting point for the analysis of the
higher genus case.  

The structure of the article is as follows. In Section~\ref{sec2} we give the 
basic definitions and background material on the phase space, dimension theory,
semi-dispersing and cylindric billiards. In Section~\ref{sec3} we collect
all the geometric arguments which can be presented in terms of the two ball
system. In particular we define the important notions of long/short symbolic collision
sequence and richness. In Section \ref{sec4} the lift is defined and analyzed.
Furthermore this section includes all the dynamical/topological arguments
and finally it contains the proof of the main theorem.

\section{Basic definitions}\label{sec2}
\setcounter{equation}{0}

\subsection{Phase space} \label{sec2.1}
A billiard is a dynamical system describing the
motion of a point particle in a connected, compact domain
$Q\subset\R^d$
or $Q\subset\T^d=\R^d/ {\mathcal L}$, $d\ge 2$, with a piecewise
$C^2$-smooth boundary. Here ${\mathcal L}\subset \R^d$ is a lattice, 
i.e. a discrete subgroup of the additive group $\R^d$ with 
${\rm rank}({\mathcal L})=d$. For the billiards considered in this paper 
${\mathcal L}$ will always be a product of several copies of two lattices: 
the hexagonal lattice in $\R^2$ and $\Z$ in $\R$.   

Inside $Q$ the motion is uniform, whereas the
reflection at the boundary $\partial Q$ is elastic (the angle of
incidence equals the angle of reflection). When a ball reaches a corner 
the collision law may
not be well defined, finitely many outgoing velocities and, correspondingly, 
finitely many trajectory branches appear (see Subsection~\ref{sec2.3} on multiple 
reflection points, i.e.\ on corner points).
Since the absolute value
of the velocity is a first integral of the motion, the phase space of a 
billiard can be identified with the unit tangent bundle over $Q$.
Namely, $M= Q\times \S^{d-1}$,
where $\S^{d-1}$ is the unit $d-1$-sphere.
In other words, every phase point $x$ is of the form
$(q,v)$ where $q\in Q$ and $v\in\S^{d-1}$ is a tangent vector
at the footpoint $q$.  The billiard flow will be denoted by $\Phi_t: M \to M$.
There is a  natural invariant probability  measure, called the Liouville measure, 
which we denote by $\mu$. It satisfies $d\mu = {\rm const.} \, dq \, dv$.

A special case is the billiard
system of two balls of equal masses in a polygon $P$. Let
$q := (q_1,q_2) \in P \times P$ be the position of the two balls,
$v := (v_1,v_2) \in \S^3 
( = \R^2 \times \R^2$ such that $|v_1|^2 + |v_2|^2 = 1)$
be their velocities, and $R=2r$ is twice the (common)
radius of the two balls.  Then the phase
space $M$ is $\big \{ (q,v): q \in P^2, |q_1 -q_2| \ge R, v \in \S^3
\big\}$,
with velocities just before and after collisions identified.

If we play billiards in a rational polygon with a single point particle,
then there is a well known construction
of a flat surface (with singularities). In the cases considered in
this article the flat surface is the torus $\T^2$, this fact is
fundamental in the proof of our theorem.
Suppose the angles of $P$ are $m_i/n_i \, 180^\circ$ where $m_i$ and $n_i$ are
coprime integers.
Let $G = G(P)$ be the dihedral group generated by the reflections in the
lines through the origin which meet at angles $\pi/N$ where $N$ is
the least common multiple of the denominators $n_i$.  We consider the phase
space $X = \{(q_1,v_1): q_1 \in P, v_1 \in \S^1\}$
of this billiard flow and let $X_{\theta}$ be the subset of points whose
velocity belongs to the orbit of $\theta$ under $G$. The set $X_\theta$
can be thought of as a flat surface (in our case the torus $\T^2$, 
thus without singularities) by gluing 
the sides
of the $2N$ copies of $P$ according to the action of $G$. For further 
details see \cite{MT} and Subsection \ref{ss1}.

In the case of two balls the group $G \times G$ 
contains the information
about the effect of the ball-to-wall collisions on the velocities,
a reflection or rotation.
More precisely let $t_i$ be the
moments of ball-to ball collisions on a trajectory 
segment and furthermore, let
$t_{i-}$ and $t_{i+}$ be non-collision  time moments directly just
before and after the given ball-to-ball collisions, respectively. Then
\be\label{GG}
v_k(t_{(i+1)-})=g_i^{(k)} v_k(t_{i+}) \hbox{ for } k=1,2
\ee 
for some $(g_i^{(1)},g_i^{(2)})\in G \times G$,
depending on the ball-to-wall collisions the two balls have before
meeting again.

\subsection{Dimension theory}\label{sec2.2}

We recall some notions from topological dimension theory, to be used later in 
Section~\ref{sec4}. For a broader exposition see \cite{ball} and references 
therein. Different notions of topological dimension coincide for 
separable metric spaces and, in particular, for compact 
differentiable manifolds 
(\cite{ball,nandbox}). Let $\dim A$ be the notation for any of these for
$A\subset \mathbf M$ where $\mathbf M$ is a compact differentiable manifold. 
Actually, 
this article uses only the concept of one and two codimensional sets, always 
characterized by means of the two lemmas below.

\begin{lemma}\label{Lemma2.10} 
For any subset $A\subset\mathbf M$ the condition $\dim A\le\dim\mathbf M-1$ is
equivalent to ${\rm int}A=\emptyset$. (See Property 3.3 of \cite{ball} or 
Lemma 2.10 of \cite{nandbox}.)
\end{lemma}

\begin{lemma}\label{Lemma2.9}
For any closed subset $A\subset\mathbf M$ the following three
conditions are equivalent:

\begin{enumerate}
\item $\dim A\le\dim\mathbf M-2$;
\item{} ${\rm int}A=\emptyset$ and for every open and connected set 
$G\subset\mathbf M$ the difference set $G\setminus A$ is also connected;
\item{} ${\rm int}A=\emptyset$ and for every point $x\in\mathbf M$
and for any neighborhood $V$ of $x$ in $\mathbf M$ there exists a smaller
neighborhood $W\subset V$ of the point $x$ such that, for every pair of 
points
$y,\, z\in W\setminus A$,
there is a continuous curve $\gamma$ in the set $V\setminus A$ connecting 
the
points $y$ and $z$.
\end{enumerate}
(cf. Property 3.4 in \cite{ball} or Lemma 2.9 in \cite{nandbox}.)
\end{lemma}

A useful notion of ``small'' set in the dynamical context is 
the following, mainly because of the property just below the definition.

\begin{definition}
\label{slim}
A subset $A$ of $\mathbf M$ is called {\bf slim} if $A$ can be covered by a
countable family of codimension-two (i.e.\ at least two), closed sets of
$\mu$-measure zero, where $\mu$ is some smooth measure on $\mathbf M$.
\end{definition}

\begin{lemma}
\label{cslim}
If $\mathbf M$ is connected, then the complement $\mathbf M\setminus A$
of a slim set $A\subset\mathbf M$ necessarily contains an arcwise connected,
$G_\delta$ set of full measure.
(cf. Property 3.6 in \cite{ball}.)
\end{lemma}

We state furthermore two Lemmas on the ``additivity of codimensions''.

\begin{lemma}
\label{adddim}
Suppose $M_1$ is a one-codimensional submanifold of $M$ and $A \subset M_1$ is closed and
has empty interior in $M_1$. Then $A$ has codimension 2.  
\end{lemma}

\begin{proof}
This Lemma is a special case of Lemma 3.13 from \cite{inv} part I.
\end{proof}

\begin{lemma}
\label{intdim}
If $A\subset M_1\times M_2$ is a closed
subset of the product of two manifolds, and for every $x\in M_1$ the set
$$
A_x=\left\{y\in M_2:\, (x,y)\in A\right\}
$$
is slim in $M_2$, than $A$ is slim in $M_1\times M_2$.
(cf.\ Property 3.7 from \cite{ball}.)
\end{lemma}

\subsection{Semi-dispersing billiards} \label{sec2.3}  A billiard 
(on billiards in general, see subsection~\ref{sec2.1}) is called semi-dispersing 
if any smooth component of 
the boundary $\partial Q$ is convex as seen from the outside of 
$Q$. Semi-dispersing billiards are  typically hyperbolic systems with singularities.
We only give a short discussion of these phenomena, for a 
detailed exposition see \cite{BCST} or \cite{KSSz}. All our arguments on semi-dispersing
billiards are self contained.

There are two possible types of {\bf singularities} for billiards. 
A collision at the boundary point $(q,v)\in \partial M$ is said to 
be {\it
multiple} if at least two smooth pieces of the boundary $\partial Q$ 
meet at
$q$, and is {\it tangential} if the velocity $v$ is tangential to $\partial
Q$ at $q$.
At tangential reflection points the flow is continuous, though not
smooth, while at multiple reflection points it is not even
continuous. Thus the future semitrajectory (or the outgoing velocity) is not
well-defined for a multiple reflection point -- for such points two {\it
trajectory branches} can be considered as the limits of the smooth
dynamics. We shall denote the set of all singular reflection points
(belonging to any of the above two types, in case of multiple collision
supplied with outgoing velocity $v^+$) by $S^+$.

We introduce some more notation. Phase points with at most one singular 
reflection on their entire trajectory will be referred to as $M^*$. 
Those for which the entire orbit avoids $S^+$ are often called {\bf regular}.
This set is denoted by $M^0$ while $M^1=M^*\setminus M^0$ contains points 
with exactly one singular reflection on their orbit.
We recall the following crucial facts:

\begin{lemma}
\label{singslim}
The set $M\setminus M^0$ is a countable union of manifolds of codimension at 
least one, thus it has zero $\mu$-measure. $M\setminus M^*$ is a 
countable union of manifolds of codimension at 
least two, thus it is slim. (cf. Lemma 2.11 from \cite{nandbox}.)  
\end{lemma} 

The treatment of {\bf hyperbolicity} in semi-dispersing billiards is 
traditionally related 
to {\it local orthogonal manifolds} (or {\it
fronts}) and {\it sufficient phase points}. 

Let $x=(q,v)\in M\setminus \partial M$ and consider a $C^2$-smooth codimension
1 submanifold $W'\subset Q\setminus \partial  Q$ such 
that $q\in W'$
and $v=v(q)$ is the normal vector to $W'$ at $q$. We define 
$W$, a section of the unit tangent bundle on $Q$ restricted to
$W'$, by picking the unit normal vector for any point of $W'$. 
The section $W$ is called a {\bf local orthogonal manifold} or simply 
a {\bf front}. A
front is said to be (strictly) convex whenever its second fundamental 
form $B_{W'}(y)$ 
is positive semi-definite (positive definite) in every point $y\in W'$. 
The definition of (strictly) concave fronts is analogous.  

To arrive at sufficient phase points we first define the neutral subspace for 
a {\it non-singular} trajectory segment $\Phi^{[a,b]}x$.
Suppose that $a$ and $b$ are {\it not moments of collision}.

We will call a point $x \in M$ {\bf hyperbolic} if it has exactly one zero 
Lyapunov exponent (i.e.~the flow direction).  For almost all
hyperbolic points unique
local stable (unstable) manifolds of positive inner radius exist, these
are strictly concave (convex) fronts \cite{KSSz}.

\begin{definition} 
The {\bf neutral space} $\mathcal N_0(\Phi^{[a,b]}x)$
of the trajectory segment $\Phi^{[a,b]}x$ at time zero ($a<0<b$) is
defined by the following formula:

\begin{eqnarray*}
{\mathcal N}_0 (\Phi^{[a,b]}x):=&\{& \hspace{-0.3cm} w\in \R^d : \exists(\delta>0) s.t. \forall
\alpha \in  (-\delta,\delta)\\ 
& & \hspace{-0.3cm} v\Big(\Phi^a\big(q(x)+\alpha w, v(x)\big) \Big)=v(\Phi^ax) \ \& \\ 
& & \hspace{-0.3cm} v\Big(\Phi^b\big(q(x)+\alpha w, v(x)\big)\Big)=v(\Phi^bx) \}.
\end{eqnarray*}
\end{definition}

Observe that $v(x) \in {\mathcal N}_0 (\Phi^{[a,b]}x)$ is always true, the
neutral subspace is at least $1$ dimensional. Neutral subspaces at time
moments different from $0$ are defined by ${\mathcal N}_t (\Phi^{[a,b]}x):={\mathcal N}_0 (\Phi^{[a-t,b-t]}(\Phi^tx))$, thus they are naturally isomorphic to the one at $0$. Having the trajectory segment fixed we often use this isomorphism to omit subscripts and refer to the neutral subspace simply as $\mathcal N$. 

\begin{definition}
\label{suff}
The non-singular trajectory segment $\Phi^{[a,b]}x$ is {\bf
sufficient} if for some (and thus for any) $t \in [a,b]$: 
${\rm dim}({\mathcal N}_t (\Phi^{[a,b]}x))=1$. A regular phase point $x$ is 
said to be sufficient
if its entire trajectory $\Phi^{(-\infty,\infty)}x$ contains a finite 
sufficient segment. 
\end{definition}

Singular points are treated by the help of trajectory branches (see above): 
a point $x\in M^1$ is {\em sufficient} 
if both of its trajectory
branches are sufficient.

Sufficiency has a picturesque meaning; roughly speaking a trajectory segment
is sufficient if it has encountered all degrees of freedom. Nevertheless the
concept is important as very strong theorems hold in open neighborhoods of 
sufficient points.

\begin{theorem}{\rm ({\bf Local Hyperbolicity Theorem}, \cite{SC}.)}
\label{lochyp} Every
sufficient phase point $x \in M^0$ has an open neighborhood $x\in U\subset M$, 
such
that $\mu$ a.e. $y\in U$ is hyperbolic. 
\end{theorem}

The even more important local ergodicity theorem dates back to the articles
\cite{SC} and \cite{KSSz}. For a detailed discussion we refer to \cite{BCST}. 
Below two conditions are given under which the theorem can be proved as shown 
in \cite{BCST}. 

We need to fix some terminology first. The zero-set of a system of 
polynomial
equations in $\R^n$ is an {\it algebraic variety} (we will use these notions 
over the real ground field). Any (measurable) subset of an algebraic variety
will be called an {\it algebraic subset}. Dimension (codimension) of an
algebraic variety (and, correspondingly, of an algebraic subset) is understood
in the following sense. Consider the ideal of polynomials vanishing on the
variety and a minimal number of polynomials $P_1,...P_r$ generating that
ideal. Dimension is the maximum (taken over all points of the variety) of
$n-m$ where $m$ is the rank of the matrix $[\partial P_1,..., \partial P_r]$,
calculated at any point. We use this notion of dimension
only to formulate the condition below. 

\begin{condition}
\label{algebra}
{\rm The semi-dispersing billiard is algebraic} in the sense that
 $\partial Q$ is a finite union of one-codimensional algebraic
subsets (as subsets of $\T^d\subset \R^d$).
\end{condition}

For the second condition one more notation is introduced. Let us denote by 
$m_{S^+}$ the induced Riemannian measure on the set
of singular reflections $S^+$.

\begin{condition}
\label{scansatz}
{\rm ({\bf Chernov-Sinai Ansatz}, cf.\ Condition 3.1 from \cite{KSSz}.)}
For $m_{S^+}$-almost every point $x \in S^+$ we have $x\in
M^*$ and, moreover, the positive semitrajectory of the point $x$ is
sufficient. 
\end{condition}

The following local ergodicity theorem is the combination of
three theorems, Theorem 5.13, Theorem 5.7 and Theorem 4.4, in
\cite{BCST}.
\begin{theorem}{\rm ({\bf Local Ergodicity Theorem or Fundamental Theorem of
Se\-mi-Dis\-per\-sing Billiards}.)}
Consider a se\-mi-dis\-per\-sing billiard which is {\it algebraic} and satisfies
the {\it Chernov-Sinai Ansatz}.
Then every sufficient phase point $x \in M^*$ has an open neighborhood  
$x\in U\subset M$ which belongs to one ergodic component. 
\end{theorem}

\begin{remark}
Recent research by Sim\'anyi \cite{nlocerg} indicates that the condition of 
algebraicity may 
not be necessary for the local ergodicity theorem. Nevertheless, billiards 
discussed in this article are all algebraic, thus the version above is 
applicable.
\end{remark}

\medskip

{\bf Cylindric billiards} make an important subclass of semi-dispersing ones. 
In their setting the
configuration space is defined by cutting out a finite number of 
cylindric regions from
the $d$-dimensional unit torus, i.e.
$Q=\T^d \setminus (C_1 \cup\cdots\cup C_k)$ where $\T^d=\R^d \big/ {\mathcal L}$. 
For the precise definition of the cylinders we need
three data for each $C_i$. We fix $A_i$, a subspace of the $d$-dimensional 
Euclidean
space $\R^d$, the {\bf generator subspace} of the
cylinder. The subspace $A_i$ should be a
lattice-subspace  (i.e.\ the discrete intersection $A_i\cap {\mathcal L}$
has rank equal to $\dim A_i$, cf.\ \cite{trans,ncyl})
to get a properly defined cylinder on
$\T^d$ after factorization.
We assume $\dim (L_i)\ge 2$, where $L_i=A_i^{\perp}$ is the notation for the 
{\bf base subspace}, the orthogonal complement of $A_i$. 
The {\it base}, $B_i \subset L_i$ is a convex, compact domain, for
which the $C^2$-smooth boundary $\partial
B_i$ is assumed: (i) to have everywhere positive definite second fundamental
form (to ensure semi-dispersivity),
and (ii) to be a one-codimensional algebraic subset of $L_i$ (to ensure
algebraicity of the billiard, i.e.\ the validity of
Condition~\ref{algebra}). Furthermore a translational vector $t_i \in \R^d$
is given to place our cylinder in  
$\T^d$. By the help of these data our cylinders are defined as: 
\be
C_i:=\{a+l+t_i : a\in A_i, l\in B_i \} \big/ {\mathcal L}.
\label{defcyl}
\ee

The most important conjecture related to cylindric billiards is the following 
one. 
\begin{remark}{\rm ({\bf Transitivity conjecture}.)} \label{trans}
The cylindric billiard is hyperbolic and ergodic if and
only if the system of base subspaces $L_1, \cdots ,L_k$ has the Orthogonal
Non-Splitting Property. That is there is no orthogonal splitting $\R^d=K_1
\oplus K_2$ for which $dim(K_j)>0 \, (j=1,2)$ and which has the property that 
for any
$i=1, \cdots ,k$ either $L_i \subset K_1$ or $L_i \subset K_2$. 
\end{remark}
The name for the conjecture is related to the fact that the condition of 
orthogonal non-splitting can be equivalently stated in terms of the 
transitivity  of a  certain group action. Strongest results related to it were 
published \cite{ncyl}. More details on cylindric billiards can be found in 
\cite{ncyl,trans,bcyl}.  

\medskip

Before closing the section we recall two theorems from the survey \cite{ball}. For 
both of them we are given a flow $(X,\phi^t)$ with invariant measure $m$ and a set 
of time moments $H$. Given $B\subset X$ the notation $A_H(B)$ refers to the set
$$
\{ \ x\in X \ | \ \phi^t(x)\cap B = \emptyset \ \forall t \in H \ \}. 
$$
    
\begin{theorem}\label{ballav}
{\rm ({\bf Weak ball avoiding theorem}.)}
Assume that the flow $\phi^t$ is mixing, $\sup H = +\infty$ and $m(B)>0$. 
Then $m(A_H(B))=0$.
\end{theorem}

Theorem~\ref{ballav} holds for any mixing flow, not necessarily a billiard flow.
In the second theorem we restrict ourselves to the case of semi-dispersing 
billiards. For a more general formulation we refer to \cite{ball}.  

\begin{theorem}{\rm ({\bf Strong ball avoiding}.)}
\label{sballav} 
Assume that $X$ is a full-measure invariant set in a semi-dispersing 
billiard that (i) satisfies the conditions of
the local ergodicity theorem and for which 
(ii) the complement of sufficient points is slim. Furthermore,
$B\ne\emptyset$ is open and for the set of time moments $H$ we have 
$\inf H =-\infty$, $\sup H=\infty$. Then $A_H(B)$ is slim. 
\end{theorem}

\begin{remark}
\label{strongball}
If conditions (i) and (ii) above hold for all points of the phase space the
semi-dispersing billiard is automatically mixing. 
On the other hand, in certain orthogonal cylindric billiards 
known to be mixing some trivial one-codimensional submanifolds appear 
consisting of trajectories that do not collide with all the cylinders, 
and are, consequently, non-sufficient (see also Remark~\ref{square} and 
references \cite{ball,S,trabant}). In such a case it is the complement of 
these trivial trajectories -- an invariant set of full measure -- to which 
Theorem~\ref{sballav} applies. 
\end{remark}

\section{Downstairs}\label{sec3}
\setcounter{equation}{0}

Recall the notion of neutral subspace from subsection~\ref{sec2.3}.
We give one more definition, that of the {\it advance}. Consider a
non-singular orbit segment $\Phi^{[a,b]}x$ with a collision
$\sigma$ taking place at time $\tau=\tau(x,\, \sigma)$. For
$x=(q,v)\in M$ and $w\in\R^4$, $\Vert w\Vert$ sufficiently small,
introduce the notation $T_w(q,v):=(q+w,v)$.

\medskip

\begin{definition} 
For any collision $\sigma$ of $\Phi^{[a,b]}x$
and for any $t\in[a,b]$, the {\bf advance}
$$ \alpha_\sigma :\, \mathcal N_t(\Phi^{[a,b]}x)\rightarrow\R
$$ is the unique linear (in $w$) functional which  satisfies 
$$ \alpha_\sigma (w):=\tau(x,\sigma)-\tau(\Phi^{-t}T_w \Phi^tx,\sigma)
$$
in
a sufficiently small neighborhood of the origin of $\mathcal
N_t(\Phi^{[a,b]}x)$.

\end{definition}

\medskip

With the help of the advance we may have a more explicit description
of  neutral vectors.

Namely, consider any ball-to-ball collision $\sigma$ of the trajectory
$\Phi^{[a,b]}x$, any fixed time moment $t$ (close enough to
$\tau(x,\sigma)$), and any neutral vector $w\in\mathcal
N_t(\Phi^{[a,b]}x)$ (with $||w||$ small enough). For more explicit
notation $v=(v_1,v_2)$ and $w=(w_1,w_2)$, where the subscripts
indicate the two-dimensional components of the velocity and the
neutral vector, corresponding to  the first and the second ball,
respectively. All of them are considered in the time moment $t$. There 
are two different kinds of neutral perturbations, the first is 
translating each of the two balls by the same vector $n \in \R^2$, 
while the second is moving along the flow direction by the advance. 
This yields:

\be
\label{neutral}
w=(w_1,w_2)=(n,n)+\alpha_\sigma (w) \cdot (v_1,v_2).
\ee

\medskip

\begin{convention}
\label{goodsegment}
From this point on throughout the article we consider trajectory segments $\Phi^{[a,b]}x$ 
for regular phase points
$x$ for which the first and last collisions are both ball-to-ball.
\end{convention}

Let us fix a segment $\Phi^{[a,b]}x$ of the above type. 
Time moments for the consecutive ball-to-ball collisions 
will be denoted by $t_1,\dots,t_{k+1}$. Otherwise we use
the notations of Equation \eqref{GG}.
The {\bf long symbolic collision sequence} corresponding to $\Phi^{[a,b]}x$ is
$$b \g_1 b \g_2 \dots \g_k b$$
where $\g_i = (g_i^{(1)},g_i^{(2)})$ and $b$ denotes the ball-to-ball collisions.
We define $\g$ to be {\bf simple} if $g := g^{(1)}=g^{(2)}$.

For brevity in the arguments below we use the notation $\alpha_i$ for the 
advance of the collision at time moment $t_i$.
Consider the particular case of two consecutive ball-to-ball collisions at $t_j$ and $t_{j+1}$ 
with $\g_j$ simple. Given any
neutral vector formula (\ref{neutral}) applies for the 
collision at time moment $t_j$ (more precisely, for non-collision times just 
before or after the collision).
\be
\label{neutral1}
w=(w_1,w_2)=(n,n)+\alpha_j (v_1,v_2),
\ee
and similarly for the collision at $t_{j+1}$:
\be
\label{neutral2}
w'=(w_1',w_2')=(n',n')+\alpha_{j+1} (v_1',v_2').
\ee
By dynamics 
$$ 
(v_1',v_2')=(gv_1, gv_2), \qquad (w_1',w_2')=(gw_1, gw_2)
$$ 
thus necessarily $n'=gn$ and $\alpha_j=\alpha_{j+1}$. To see this note that 
$v_1'=v_2'$ is not possible right before or after a ball-to-ball collision. 

Now we will begin to define {\em short collision sequences}.
In view of the above, if we have $b \g_1 b \dots \g_k b$ with all the $\g_i$ simple,
then dynamically this has the same effect as a single ball-to-ball collision. The
only role of the $\g_i$ that we need to note is that a neutral vector of the form
$(n,n)$ just before the first ball-to-ball collision evolves into $(n',n')$. Here 
$n' = sn$ with $s = \prod_{i=1}^k g_i^{(1)} =  \prod_{i=1}^k g_i^{(2)}$.
In accordance with
this fact we shall use the symbol $(b,s)$ for a maximal sequence of
ball-to-ball collisions such that the consecutive ones are all
separated by simple $\g$-s. In case $k=0$, i.e.\ if the sequence consists of a 
single ball-to-ball collision we fix $s=Id$. Following tradition (e.g. \cite{inv}) 
such a maximal sequence will be referred to as an {\bf island}. As discussed above, 
ball-to-ball collisions in an island have the same advance, thus we may define {\it the unique
advance for the island}. 

Consider a trajectory segment whose long symbolic collision sequence satisfies 
that neither $\g_1$ nor $\g_k$ are simple. 
In view of all the observations made above we use the notation
$$b\,\g_1\,(b,s_1)\,\g_2\,(b,s_2)\dots (b,s_{K-1})\,\g_K\,b$$ 
for {\bf the short symbolic collision sequence} of the trajectory 
segment. Here the symbol $(b,s_i)$ refers to an island (see above) while the 
(non-simple pairs of) group elements $\g_i$ describe the effect of the ball-to-wall collisions in
between. There is a slight ambiguity of notation, the $\g_i$ in the short symbolic
collision sequence is a subsequence of the $\g_j$ for the long one. This ambiguity
should not cause any confusion.  Note that the short symbolic collision sequence consists of
$K+2$ islands since the first and last $b$, which each denote a single ball-to-ball 
collision, are also considered as islands.


We define two more symbols.  
Let $\hg_i$ and $\bg_i$  be the unique elements of $G$ for which
$g_i^{(2)}=\hg_ig_i^{(1)}$ and $g_i^{(2)}\bg_i=g_i^{(1)}$. 
The transformation $\hg_i$ tells us how the relation
of the two velocity vectors has changed between the two ball-to-ball
collisions. We have $(g_i^{(2)})^{-1} = \bg_i (g_i^{(1)})^{-1}$, thus
the transformation $\bg_i$ plays the same role for the backwards dynamics.

\begin{remark}
\label{timeref}
Symbolic collision sequences defined this way are not time reflection 
symmetric.
We note that for the case of the square, as $G$ is commutative, $s_i^{-1}=s_i$ and 
$\hg_i=\bg_i = (\bg_i)^{-1}$ 
automatically. This fact enabled Sim\'anyi
to use another concept of symbolic collision sequence, which is time reversal symmetric, 
in his article \cite{nandbox}. 
\end{remark}

Any group element $g \in G$ may be either a reflection (in this
case we use the notation $g=R$) or a rotation ($g=O$). Sometimes it is useful 
to indicate $R_E$ where the line $E\in \R^2$ is the axis of the reflection 
$R$. 
If $\hg_i$ is a reflection/rotation then the same is true for $\bg_i$.

\begin{definition}\label{rich}
A long symbolic collision sequence is {\bf rich} if it contains a subsequence
with short form  
$b\,\g_1\,(b,s)\,\g_2\,b$
where either
\begin{enumerate}
\item $\hg_1=R_E$ and $\bg_2=R_{E'}$ with $E' \ne sE$, 
or
\item $\hg_1=R$ and $\bg_2=O$ for any reflection and rotation.
\end{enumerate}
A phase point $x\in M^0$ is rich if its entire trajectory contains a finite segment with 
rich collision sequence. For rich points of $M^1$ the same should hold for both trajectory 
branches.  
\end{definition}

\begin{lemma}
\label{hypgeom}
For any phase point $x\in M$ with rich collision sequence there 
exists a neighborhood $U$ and a one-codimensional submanifold
$L\subset M$ such that any $y\in U\cap (M \setminus L)$ is 
sufficient.
\end{lemma}

Let us describe the neutral vectors for long collision sequences of length three 
first.

\begin{sublemma}
\label{Sref}
Assume $x$ has a trajectory segment with long collision sequence of the form 
$b\,\g\,b$  such that $\hg=R$ for some reflection $R=R_E$. Consider a neighborhood 
$U$ of $x$, such that for $y \in U$  
(this finite segment of) the collision sequence is the same.
Then, apart from a degeneracy (present on a one-codimensional
manifold $L \subset U$) for any $y \in U \backslash L$ the neutral 
subspace $\mathcal N$
is two dimensional
and for any vector in $\mathcal N$
the advances of the two ball-to-ball collisions are equal to one another.
\end{sublemma}

\begin{proof}
Let us fix non-collision time moments just after the first and just before
the second ball-to-ball collisions and denote them with  
$t^*$ and $t^-$, respectively. At these 
time
moments any vector of $\mathcal N$, by neutrality with respect to the
ball-to ball collisions, has the form: 
\bea
\label{twocoll}
w^*=(w^*_1,w^*_2) &=& (m,m)+ \alpha (v^*_1,v^*_2), \\
\label{twocoll2}
w^-=(w^-_1,w^-_2)
&=& (n,n)+ \beta (v^-_1,v^-_2).  
\eea
Here $\alpha$ and $\beta$ are
the advances of the two ball-to-ball collisions, the upper indices
indicate time moments for the two dimensional velocity vectors/neutral
vectors of the two balls, while $m,n\in \R^2$ are arbitrary, cf. 
(\ref{neutral}).
We need
to evolve the first of these two equations to the time moment $t^-$ to
compare it with the second: $v^*$ turns into $v^-$, while, as $\hg=R$, 
$(m,m)$
evolves into $(l,Rl)$ for some $l\in\R^2$. Thus we have:  
\be\label{3.3}
w^-=(w^-_1,w^-_2) = (l,Rl)+ \alpha (v^-_1,v^-_2).   
\ee
There are
two possibilities. Either $\alpha\ne\beta$. Comparing 
Equations~\eqref{3.3} and
\eqref{twocoll2} implies
\be
l-Rl=(\alpha-\beta)(v^-_1-v^-_2) 
\ee
that results in $(v^-_1-v^-_2)\in
E^{\perp}$ where $E^{\perp}$ is the line in 
$\R^2$ perpendicular to $E$, to the axis of the
reflection $R$. This means $x\in L$ where $L$ is a one-codimensional
submanifold of $M$.  On the other hand if $\alpha=\beta$, again comparing
Equations~\eqref{3.3} and \eqref{twocoll2}
yields $n=l=Rl$, thus the a priori two  dimensional $n$ is restricted to $E$, 
to the axis of $R$.
\end{proof}

\begin{sublemma}
\label{Srot}
Assume $x$ has a trajectory segment with long collision sequence of the form 
$b\,\g\,b$  such that $\hg=O$ for some rotation $O$. 
Then, if the two advances are equal, the phase point $x$ is
sufficient. If the two advances are not equal, they give a full
description of the neutral subspace.
\end{sublemma}

\begin{proof}
With notations and argumentation of the previous sublemma we have the
validity of (\ref{twocoll2}) and
\be  
\label{rotcoll}
w^-=(w^-_1,w^-_2) = (l,Ol)+
\alpha (v^-_1,v^-_2).
\ee  
If $\alpha=\beta$, we have $l=n=Ol$.
However, rotations have no fixed points (except for the origin), thus 
$l$ is zero and the
neutral vector is trivial: the phase point is sufficient. If
$\alpha\ne\beta$, we get  \be  l-Ol=(\alpha-\beta)(v^-_1-v^-_2).   \ee
As the linear map $Id-O$ is invertible, $l$ (and consequently the
neutral vector) is completely determined by the advances (and the
velocity components which, however, do depend only on the phase point
itself and not on the perturbation).

{\it Alternatively}, for future reference in case of unequal advances we 
may 
derive from (\ref{twocoll2}) and (\ref{rotcoll})
\be
\label{rotdeg}
n-O^{-1}n = (\alpha-\beta) (v_1^--O^{-1}v_2^-).
\ee
With the reasoning given above $n$  is completely determined by the right 
hand side, i.e.\ the advances determine the perturbation. 
\end{proof}

\begin{proofof}{Lemma~\ref{hypgeom}}

\noindent
We must analyze short collision sequences of length five, 
$b \, \g_1 \, (b,s) \, \g_2 \, b$.

(1) Let us consider case (1) from Definition~\ref{rich} first. We denote the advances
    of the three islands as $\alpha$, $\beta$ and
    $\gamma$ and fix time moments $t_-$ and $t_+$ just before and
    after the middle island. We may apply
    Sublemma~\ref{Sref} for both trajectory segments (up to the first and
    from the last ball-to-ball collision of the island $(b,s)$) to conclude that apart 
    from codimension one
    $\alpha=\beta=\gamma$. Any neutral vector is of the form 
    \be
    \label{before}
    w^-=(w^-_1,w^-_2) = (n,n)+ \beta (v^-_1,v^-_2)
    \ee 
    at time $t_-$,
    with $n\in E$. This neutral vector evolves into 
    \be
    \label{after}
     w^+=(w^+_1,w^+_2) = (n',n')+ \beta (v^+_1,v^+_2)
    \ee 
    where $n'=sn$. 
    Applying Sublemma~\ref{Sref} in {\it backward time} to the second segment yields 
    $n'\in E'$ by
    $\bg_2=R_{E'}$. As we assumed $E'\ne sE$ and 
    we have $n\in E$, this means $n=0$ which is equivalent to 
    sufficiency in view of (\ref{before}).

(2) Now consider case (2) of Definition~\ref{rich}. For the time moments 
    and the advances we use the notations of (1). We apply 
    Sublemma~\ref{Sref} to the segment ending with the first ball-to-ball 
    collision of the island $(b,s)$ and apply,
   with time direction reversed, Sublemma~\ref{Srot} to
    the segment starting with the last 
    ball-to-ball collision of the middle island.
    Formulas 
    (\ref{before}) and (\ref{after}) are both valid. 
    Furthermore, from Sublemma~\ref{Sref} we know that (apart from
    codimension 1) $n$ is restricted to a line. 
    Thus $n'=sn$ is restricted to a line as well. 
    On the other hand, by Sublemma~\ref{Srot}, we
    may assume $\beta\ne\gamma$ (otherwise we have
    sufficiency), thus the following formula (the time-reversal 
analogue of (\ref{rotdeg})) is valid  for $v^+$: 
\be
\label{rotdegbar}
n'-O^{-1}n' = (\beta-\gamma) (v_1^+-O^{-1}v_2^+).
\ee

The left hand side of (\ref{rotdegbar}) is the image of $n'$ by the
    invertible linear map $Id - O^{-1}$. Thus it also is restricted to a line,
    inheriting this property from $n'$.
    This however means $v_1^+-O^{-1}v_2^+$ is
    restricted to a line, which implies that the phase point $x$ belongs
    to a one codimensional submanifold $L$ of the phase space.

For future reference we make one more simple remark. Let us fix a non-collision
time moment $t_{\#}$ just before the third island. Then
$v_1^{\#}=g_2^{(1)}v_1^+$ and $v_2^{\#}= g_2^{(2)}v_2^+$. On the other hand, 
$\bg_2 =O$ and thus $(g_2^{(2)})^{-1}=O(g_2^{(1)})^{-1}$. 
Combining the last three formulas yields 
$$
v_1^+-O^{-1}v_2^+=(g_2^{(1)})^{-1} (v_1^{\#}-v_2^{\#}).
$$
Thus the points of the one-codimensional degeneracy submanifold $L$ are 
equivalently characterized by the relation: $v_1^{\#}-v_2^{\#}$ is restricted 
to a line. 
\end{proofof}

\begin{remark}
\label{firstcodim}
Having a look at the arguments above it is useful to note that in all cases 
the degeneracy submanifolds are characterized by the following relation: the difference of the 
velocities of the two balls, $v_1-v_2$, is restricted to a line 
when calculated at a time moment just before a given ball-to-ball collision. 
\end{remark}

\noindent
The following Lemma plays an important role in establishing that the set of
non-sufficient points is slim.

\begin{lemma}
\label{transv1}
Consider a trajectory segment with a ball-to-ball collision on it and let 
$t_+$ and $t_-$ be non-collision time-moments: $t_-$ just before the ball to 
ball collision and $t_+$ any time moment after the ball-to-ball collision. Fix 
furthermore two arbitrary lines, $E^+$ and $E^-$ in $\R^2$. 
The set of points  for which both $v_1^--v_2^-\in E^-$ and $v_1^+-v_2^+\in E^+$ 
belongs to a two-codimensional submanifold of $M$.
\end{lemma}

\begin{proof}
Phase points with any of the above two degeneracy relations form 
one-codimensional submanifolds of the phase space. Let us denote these 
with $L_-$ and $L_+$: our task is to show the transversality of these 
manifolds. We argue along the lines of the proof of Lemma 3.10 from 
\cite{nandbox}. Fix $x_0\in L_-$ and sufficiently small numbers $
\delta,\epsilon>0$ such that for all points $x$ of
$$
U_-= \{ \ x\in L_- \ |\ d(x_0,x) <\epsilon \ \}
$$
the trajectory segment $\Phi^{[0,\delta]}x$ is collision free (i.e.\ the time
moment $t_-+\delta$ is still before the ball-to-ball collision $t_-$ precedes).
We foliate $U_-$ with convex local orthogonal manifolds (i.e.\ with fronts). 
Namely, consider 
the equivalence relation defined for $x,y\in U_-$ as
\be
\label{eqclass}
x \sim y \ \Longleftrightarrow \ (q_1^-(x)-q_1^-(y)) = (q_2^-(y)-q_2^-(x)) \perp E^-.
\ee
The equivalence classes $C(x)$ of $\sim$ are 3 dimensional submanifolds of 
$U_-$ (2 dimensional in the velocity space and 1 dimensional in the 
configuration space). Furthermore, for small positive times 
$0<t<\delta$ they evolve into 
convex local orthogonal manifolds: each $\Phi^t C(x)$ is a front 
strictly convex in a two dimensional plane, the only neutral direction is 
\be
\label{neu}
(w_1,w_2); \quad w_1 = - w_2 =:w \perp E^- .
\ee
A perturbation of the form above is definitely not neutral with respect to 
the ball-to-ball collision the time moment $t_-$ precedes. To see this recall 
that 
(i) our phase point $x$ belongs to $L_-$ thus the vector (\ref{neu}) is 
perpendicular to the flow direction (i.e.\ to the velocity), (ii) any 
perturbation 
neutral with respect to the ball-to-ball collision should have 
the form (\ref{neutral}).

As a consequence $C(x)$, when considered after the ball-to-ball collision, evolves into a convex local orthogonal manifold which is strictly convex in all the three dimensions. Strict convexity of local orthogonal manifolds is preserved by the flow, thus $\Phi^{t_+-t_-}C(x)$ (the front considered at $t_+$) is strictly convex. This, however, means that it is necessarily transversal to $L_+$ which is defined by linear restriction on the velocity. 
\end{proof}

\begin{definition}
A point in $M$ is called {\bf twice rich} if its orbit contains two trajectory 
segments with long symbolic collision sequences 
that, on the one hand, may intersect in at most one symbol $b$, and on the other hand, 
\begin{enumerate}
\item are both rich in forward or both rich in backward time, or 
\item the first one is rich in forward time while the second is rich in 
backward time.
\end{enumerate} 
\end{definition}

\begin{corollary}
\label{2rich}
Those phase points that are twice rich and non-sufficient form a slim subset.
\end{corollary}

\begin{proof}
For those points which satisfy the first condition above this is a 
straightforward  consequence of Lemmas~\ref{hypgeom},~\ref{transv1} and 
Remark~\ref{firstcodim}. In the case of opposite orientations sufficiency 
follows from the second sequence (Lemma~\ref{hypgeom} applied in backward time)
unless $v_1-v_2$ is restricted to a line in a non-collision time moment $t_+$, 
where $t_+$ is after 
the last ball-to-ball collision on the first sequence. Thus 
Lemma~\ref{transv1} applies even in this case.
\end{proof}

\section{Upstairs}\label{sec4}
\setcounter{equation}{0}
\subsection{Lifting}\label{ss1}

Let $P$ be an integrable polygon and $M$ the phase space of two disks in 
$P$.
We consider the billiard flow $\Phi_t: M \to M$ which we refer to as the
polygonal flow.  We want to lift this system to a cylindric billiard
flow
$\Psi_t: N \to N$ where $N$ is the four dimensional torus $\T^4$ with some
cylindric scatterers removed. The projection $\pi:N \to M$ will be a 
continuous, measure
preserving and finite to one semi-conjugacy.  The rest of this subsection 
defines 
$\pi$ and
proves its important properties, in particular that $\pi$ ``preserves 
codimension'' (see Lemma~\ref{lemma1}).

In the next subsection we will use this ``preservation of  
codimension'' to
prove hyperbolicity and ergodicity of the billiard system in the polygon.

If we consider only one point particle in $P$, then there is a natural unfolding 
process: instead of
reflecting the ball when it collides with the boundary of $P$, reflect $P$ 
in the side of collision
and continue the orbit of the ball in a straight line.  If we fix a 
``generic''
initial direction $\theta$
of the ball in $P$, then there are $|G|$ possible directions the orbit can 
take.
Taking one copy of $P$ for each direction, and gluing together parallel 
sides via the unfolding procedure, yields the two torus 
$\T^2=\R^2 \big/ {\mathcal L}$. Here the lattice ${\mathcal L}\subset \R^2$ is 
either $\Z^2$ (if $P$ is the square or the half square) or the 
hexagonal lattice (if $P$ is the equilateral triangle or the half equilateral triangle). 
The corresponding billiard flow 
decomposes into
a one-parameter family of linear flows on $\T^2$, see \cite{MT} for details.  

This way we may think of $\T^2$ as a union of $|G|$ copies of $P$.
Accordingly, we use the notation $z = gq$ for $z \in \T^2$, where $g \in G$ and $q \in P$.
Furthermore, for $h \in G$ fixed, let $P_h = \{z=gq \in \T^2: g = h\}$.
The {\em cells} $P_h$ overlap only at their boundaries. 

Note that the motion of one hard ball in the polygon $P$ is equivalent to the motion 
of a point particle in a smaller copy of $P$, thus the above unfolding process applies.
We generalize this to two balls in $P$.  First pretend that the two balls 
are transparent, i.e.~pass
through each other.  Then the above  construction yields the linear flow 
on $\T^4 = \T^2 \times
\T^2$.  We can think of $\T^4$ as the union of cells $P_{g_1} \times P_{g_2}$ with
$(g_1,g_2) \in G \times G$.
Now we must  take into account the collisions.  Clearly the 
cylinder $C_e$ (see below), corresponding to a 
neighborhood of the diagonal, projects down to overlapping balls, and thus 
does not belong to $N$.
The further cylinders which correspond to overlapping balls can be 
obtained from $C_e$ by 
the action of the symmetry group $G$ on one of the two coordinates of 
$\T^4$.

To be more precise we use coordinates $z := (z_1,z_2)$ on $\T^4$ where 
$z_i\in\T^2$.
We have $z_i=g_iq_i$ where $g_i\in G$; $q_i\in P$. The domain to be removed 
from $\T^4$ is  $\{z: |q_1 - q_2| \le R\}$ where $R$ is twice the (common) radius 
of the balls. This set is a finite union of cylindric scatterers $C_g$, 
$g\in G$ where 
\be\label{cylinder}
C_g := \{z: |z_1 - gz_2|\le R\}.
\ee 
In particular, 
$C_e := \{z: |z_1 - z_2| \le R\}$. We remark that the scatterer $C_e$ does not 
depend on the choice of $g_1=g_2$ since the two balls are in the same copy of 
the polygon $P_{g_1}=P_{g_2}$. This corresponds to the fact that we have only 
$|G|$ different scatterers (and not $|G|^2$): it is only the relation of the 
two discrete coordinates $g_i$ that matters.  
Note that the cylinder 
$C_g$ intersects the cell $P_{g_1} \times P_{g_2}$ if and only if
$g = g_1 g_2^{-1}$.

\begin{remark}
\label{lattice}
The configuration space of the above defined
cylindric billiard $N$ is a subset of $\T^4= \R^4 \big/ {\mathcal L}$ where the lattice 
${\mathcal L}$ is either $\Z^4$ (if $P$ is the square or the half square) or the 
product of two copies of the hexagonal lattice 
(if $P$ is the equilateral triangle or the half equilateral triangle). 
Straightforward calculation shows that for all cylinders defined by 
(\ref{cylinder}) both the generator and the base subspaces are lattice subspaces 
(see also the proof of Lemma~\ref{cylgeom} on the explicit form of these subspaces). 
\end{remark}

The phase space $M$ is $P^2 \times \S^3$ with identification at the 
boundary.
A point in $p \in N$ can be given coordinates $(q_1,q_2,v_1,v_2,g_1,g_2)$ 
with
$q_1,q_2 \in P$,
$g_1,g_2 \in G$ and $v_1^2 + v_2^2=1$. Then we define $\pi(p) = 
(q_1,q_2,g_1^{-1}v_1,g_2^{-1}v_2)$.
This map is clearly well-defined for those points where the coordinate 
system is unique.
For those $p$ for which $q_1 \in \partial P$ or $q_2 \in \partial  P$
there are $(g'_1,g'_2) \ne (g_1,g_2)$ describing $p$. Suppose for concreteness 
that $g'_1 \ne g_1$.
However we have
the two projections of $p$ coincide since by the definition of the phase 
space
the points  $(q_1,g'^{-1}_1v'_1)$ and $(q_1,g^{-1}_1v_1)$ are identified. Thus 
$\pi$ is
well defined.

The direct product of the  measure $\mu$ and the discrete uniform measure on $G 
\times G$
is an invariant measure $\nu$ for the billiard flow $\Psi$.

We say that a map preserves a property, if the image of a set satisfying a 
certain
property satisfies the same property. 

\begin{lemma}\label{lemma1} The map $\pi$
\begin{enumerate}
\item{}  is continuous,
\item{} is an at most $|G|^2$ to one map,
\item{}  is a semi-conjugacy,
\item{}  preserves codimension one subsets,
\item{}  is measure preserving, and
\item{} preserves closed codimension two subsets.
\end{enumerate}
\end{lemma}

\begin{proof}
(1-3) The map $\pi$ is clearly a continuous projection which is at most 
$|G|^2$ to one.  Because of the form of the 
projection
$\pi$ is a semi-conjugacy of the billiard flows away from the collisions, 
and an elementary calculation shows that it is also a semi-conjugacy at 
the collisions.
\\
(4) This follows by combining continuity with Lemma \ref{Lemma2.10}.\\
(5) By definition $\pi_* \nu = \mu$.\\
(6) We use the characterization of codimension 2 given in Lemma 
\ref{Lemma2.9}.
Consider a closed codimension 2 subset $A \subset N$. Since $G^2$ acts on 
$N$
by isometries that clearly preserve codimensionality of sets we can 
assume,
without loss of generality, that $A$ is $G^2$ invariant. 
We claim that $\pi(A)$ is
closed. For each $\g=(g_1,g_2) \in G \times G$ let $A_{\g} = A \cap P_{g_1} 
\times
P_{g_2} \times \S^3$. Clearly $A_{\g}$ is closed and $\pi(A_{\g})$ is closed as 
well
since $\pi|A_{\g}$ is a homeomorphism. Thus $\pi(A) = \cup_{\g \in G^2} 
\pi(A_{\g})$
is closed as well.

We use the third equivalent characterization of Lemma \ref{Lemma2.9}. The property that 
${\rm int}\, \pi(A) = \emptyset$ follows from (1).  Fix $x \in M$
and a neighborhood $V$ of $x$.
If $x \not \in \pi(A)$ then the property is trivial, thus assume $x \in 
\pi(A)$.
We choose some lift $\hat{x} \in \pi^{-1}(x) \cap A$.
Choose a neighborhood $\hat{V}$ of $\hat{x}$
such that $\pi(\hat{V}) \subset V$.  Since $A$ is codimension 2 there is 
a neighborhood
$\hat{W}$ of $\hat{x}$ which satisfies the property: for any $\hat{y},\hat{z} 
\in \hat{W}\setminus A$ 
there is an arc $\hat{\gamma} \subset \hat{V} \backslash A$ connecting 
them.
For each $\g \in G^2$ let $\hat{\gamma}_{\g} = \gamma \cap P_{g_1} \times 
P_{g_2} \times 
\S^3$.
The set $\gamma_{\g} := \pi(\hat{\gamma}_{\g})$ is a continuous curve.  
Furthermore,
since $\pi$ is well defined at the boundary of $P$, the union of these curves 
$\gamma:=\cup_{\g \in 
G^2} \gamma_{\g}$
gives a well defined continuous curve which avoids $\pi(A)$ since we have 
assumed that
$A$ is $G^2$ invariant.
Let $W = \pi(\hat{W})$.  For any two points $y,z \in W\setminus \pi(A)$ we can 
find 
lifts 
$\hat{y},\hat{z}$ in $\hat{W}\setminus A$, thus the already constructed curve 
$\gamma$ connects $y$ to $z$ in 
$V \setminus \pi(A)$.
\end{proof}

Lemma~\ref{where} and Corollary~\ref{welldef} describe in more detail
how different types of collision sequences are lifted. 

\begin{lemma}\label{where}
Consider a trajectory segment $\Phi^{[a,b]} x \subset M$ 
with long symbolic collision sequence $b \k b$ where $\k = (k_1,k_2)$.
Assume that the first ball-to-ball collision is lifted to 
$P_{g_1} \times P_{g_2}$.
Then the second ball-to-ball collision is lifted to  $P_{h_1} \times P_{h_2}$ 
where
$h_1 = g_1 k_1^{-1}$ and $h_2 = g_2 k_2^{-1}.$ Thus if the first
ball-to-ball collision is lifted to a collision with cylinder $C_g$, with
$g = g_1 g_2^{-1}$, the second is lifted to $C_h$ with $h = h_1 h_2^{-1} =
g_1 k_1^{-1} k_2 g_2^{-1}$.
\end{lemma}

\begin{proof}
Let us denote the velocities in $M$ just after the first and just before the 
second ball-to-ball collision by $(v_1,v_2)$ and $(v_1',v_2')$, respectively.
A time segment free of ball-to-ball collisions is lifted to a segment of a 
linear flow in $N$, thus the two velocities are lifted to the same vector 
$(v_1^N,v_2^N)$. By the definition of the map $\pi$ we have 
$g_j v_j=v_j^N = h_j v_j'$, 
$j=1,2$. 
On the other hand $v_j'=k_j v_j$, $j=1,2$. This yields 
$h_1 = g_1 k_1^{-1}$ and $h_2 = g_2 k_2^{-1}$. The formula for the element $h$
defining the cylinder is a straightforward consequence. 
\end{proof}

\begin{definition}
\label{go}
Let $G_O \subset G$ be the subgroup generated by rotations.
\end{definition}

\begin{corollary}\label{welldef}
1) An island is lifted to consecutive collisions with the same cylinder.\\
2) For a trajectory segment with long collision sequence 
$b \, \k \, b$, if the first $b$ is lifted
to a cylinder of the form $C_g$ with $g \in G_O$, and $\bk \in G_O$ (or equivalently
$\hk \in G_O$), then the second $b$ is lifted to $C_h$ with $h \in G_O$.\\
3)  For a trajectory segment with short collision sequence 
$b \, \k \, (b,s) \, \l \, b$, where $\hk = R_E$ and $\bl = R_{sE}$, the first 
and last
ball-to-ball collisions are lifted to collisions with the same cylinder.
\end{corollary}

\begin{proof}
1) By definition the consecutive ball-to-ball collisions in an island are 
separated by  simple group elements $\k=(k_1,k_2)$, i.e.\ $k_1=k_2$. 
Thus Lemma~\ref{where} yields $g=h$ and therefore the coincidence of 
consecutive cylinders.

2) The elements of $G_O$ are equivalently characterized as those that are 
generated by an even number of reflections. If this property holds for $g$ and 
$\bk$ then, applying Lemma~\ref{where}, this should hold for $h$ as well.  

3) Let us assume that the first ball-to-ball collision is lifted to the cell 
$P_{g_1}\times P_{g_2}$ (and thus to the cylinder $C_{g_1g_2^{-1}}$) of $N$. 
Then, applying Lemma~\ref{where} twice the first and the last ball-to-ball 
collisions of the middle island are lifted to the cells 
$P_{g_1k_1^{-1}}\times P_{g_2k_2^{-1}}$  and 
$P_{(g_1\, k_1^{-1}\,s^{-1})}\times P_{(g_2\, k_2^{-1}\,s^{-1})}$, 
respectively. Applying 
Lemma~\ref{welldef} once more the last ball-to-ball collision is lifted to the 
cell $P_{h_1}\times P_{h_2}$ with $h_j=g_j k_j^{-1} s^{-1} l_j^{-1}$, $j=1,2$.
Note that by trivial computation $s^{-1}(R_{sE})s=R_E$ while
$\bl=R_{sE}$ and $\hk=R_E$ by assumption. 
Thus the last ball-to-ball collision of the sequence is lifted to
$C_h$ with $h=h_1h_2^{-1}=g_1 g_2^{-1}=g$.
\end{proof}

The following simple Lemma describes the geometry of the 
cylinders in $N$. Recall the basic notions related to cylindric billiards 
from Subsection~\ref{sec2.3}, in paricular, the definitions of the 
generator and the base subspaces.

\begin{lemma}
\label{cylgeom}
1) For any two distinct cylinders $C_g$ and $C_{g'}$ with $g,g'\in G_O$ both 
the generator and the base subspaces are transversal.\\
2) Except for the case of the square, the collection of cylinders $C_g$, 
$g\in G_O$ is transitive in the sense of Conjecture~\ref{trans}.\\
3) Consider $C_e$ and $C_R$ for any reflection $R\in G$. These two cylinders 
are orthogonal in the sense of \cite{S}, i.e. it is possible to choose an 
orthogonal frame in $\R^4$ such that both cylinders have generator subspaces 
spanned by vectors from this frame. Moreover, the two generator subspaces 
intersect in a line.\\ 
\end{lemma}

\begin{proof}
We note that the cylinders $C_g$ have the 
following generator and base subspaces:
$$
A_g= \{ \ (w,gw) \ | \ w\in\R^2 \ \}; \qquad 
L_g= \{ \ (w,-gw) \ | \ w\in\R^2 \ \}.  
$$
All the rest is straightforward calculation. To see (2) we remark that for 
the case of the square $G_O$ consists of two elements, the identity and 
the rotation by $180^\circ$. In all other cases rotations with different degrees 
are present.
\end{proof}

We close the Subsection with an important Remark on the upstairs 
cylindric billiard for the case of the square.

\begin{remark}
\label{square}
(1) For the case of the square the upstairs dynamical system $(N,\Psi^t,\nu)$ is 
an orthogonal cylindric billiard satisfying the conditions of the
main theorem of \cite{S}. Thus it is ergodic and hyperbolic. In contrast to 
hyperbolicity and mixing, ergodicity of the 
downstairs factor is an immediate consequence.\\ 
(2) In Subsection~\ref{dynamical} both hyperbolicity and the conditions of the Local Ergodicity 
Theorem (i.e.\ the Chernov-Sinai Ansatz, Condition~\ref{ansatz}) 
are deduced for all the four polygons 
along the same lines. 
However, non-sufficient points do not form a slim subset if $P$ is chosen to be 
the square (cf.\ Lemmas~\ref{avoiding} and~\ref{regslim}), 
thus the treatment of ergodicity should be 
slightly different. Namely, as a result of 
their geometry, in orthogonal cylindric billiards  typically trivial 
one-codimensional invariant sets appear consisting of trajectories that do not collide with all 
the cylinders and are, consequenly, non-sufficient 
(see also Remark~\ref{strongball}). In principle, these may separate ergodic 
components, however, standard methods (connecting the components with orbits of positive measure) 
exclude this possiblity. For details we refer to \cite{S,trabant}. This applies to the case of the 
square, nevertheless, to keep the exposition self-contained, we do not consider this issue here. 
Ergodicity is proved for the square by direct application of \cite{S}, see part (1) of the Remark. 
\end{remark}

\subsection{Applications of lifting}
\label{dynamical}

We begin by characterizing those points whose long symbolic collision sequence 
is extendable to infinite orbits.  We remark that no trajectory has infinitely 
many collisions in a finite time interval \cite{BFK,G,V}.
\begin{lemma}
1) If a semitrajectory (positive or negative) has no ball to
ball collisions, then the whole trajectory has no ball to
ball collisions.\\
2) The set of points whose trajectories have no ball-to-ball collisions
is slim (recall Defintion~\ref{slim}).
\end{lemma}

\begin{proof}
1) Without loss of generality let us suppose the positive 
semitrajectory has no ball-to-ball collisions.  
Consider all the points $x \in N$ whose projection $\pi(x)$ has no
ball-to-ball collisions on the positive semitrajectory.
The trajectory of each such $x$ avoids all cylinders in positive time.
Thus this semi-orbit can be thought of as the orbit of a linear flow
on $\T^4$ that avoids an open set in positive time.  Thus it avoids
the same open set in negative time as well, i.e.\ its whole trajectory avoids
all cylinders and downstairs there are no ball-to-ball collisions.

2) Since the group $G$ contains rotations, by statement 1) of 
Lem\-ma~\ref{cylgeom} we  can choose two cylinders (in $N$)
whose base spaces are transverse.
Consider those $x$ for which the trajectory of $\pi(x)$ has no ball to
ball collisions.  Then the trajectory of $x$ avoids all cylinders, in
particular the above mentioned transverse ones.
We consider the orthogonal projection $p_C$ of $\R^4$ onto the base
subspace of the cylinder $C$.  The key fact here is that
the orbit $\Psi^{(-\infty,\infty)}x$ is the orbit of a linear flow on $\T^4$, 
thus it avoids $C$ if and only if 
$p_C(\Psi^{(-\infty,\infty)}x)$ avoids $p_C(C)$.
This means that $p_C(v)$ is rationally dependent. Since this happens
for two cylinders $C$ with transverse base subspaces, Lemma~\ref{lemma1} 
parts (5) and (6) imply that
the set of $\pi(x)$'s with this property is slim.
\end{proof}

For the following definition we recall the notions of long/short collision 
sequences. 

\begin{definition}\label{poor}
1) Those regular orbits that have infinitely many ball-to-ball collisions will 
be
called {\bf extendable}.  By the above Lemma, points with nonextendable orbits
form a slim subset.\\
2) A point $x \in M$ is called {\bf O-poor} if its trajectory
is extendable and if any finite segment which starts and ends with a ball-to-ball collision 
has long collision sequence 
$b\g_1b \dots b \g_nb$ where $\bg_i\in G_O$. In other words, the 
transformations $\bg_i$ (or equivalently, $\hg_i$) are either identities 
or rotations. A point $x\in M$ is called {\bf O}$^+${\bf-poor} if there exists 
$t_0$ such that, instead of the entire orbit,
the same holds for the
semitrajectory $\Psi^tx$ for $t \ge t_0$.\\
3)   A point $x \in M$  is called {\bf R-poor} if its trajectory is extendable, not O-poor 
and if, furthermore, given  any finite trajectory
segment of the entire orbit of $x$ that has short collision sequence 
$b\,\g_1 \, (b,s)\, \g_2 \, b$, there exists a reflection $R_E\in G$ such 
that  $\hg_1=R_E$ and $\bg_2=R_{sE}$. In other words, only reflections are allowed and there is no 
segment which is rich in the sense of case (1), Definition~\ref{rich}.  
The definition and
convention of $x \in M$ 
being {\bf R}$^+${\bf-poor} is analogous to O$^+$-poor.\\
We remark that the set of points satisfying any of the above notions is a closed 
and invariant set.
\end{definition}


\begin{lemma}
\label{poorup}
(1) O-poor trajectories are lifted to trajectories that collide only with the 
cylinders $C_g$, $g\in G_O$. For the \op ones the same holds for the 
appropriate semitrajectory.\\
(2) R-poor trajectories are lifted to trajectories that collide only with two 
cylinders. One of these is always $C_e$ while the other depends on the 
trajectory, nonetheless it is $C_R$ for 
some reflection $R\in G$.    
For the \rp ones the same holds for the 
appropriate semitrajectory.\\  
\end{lemma}

\begin{proof}
In all cases, fix one particular ball-to-ball collision of the
(semi)tra\-jectory. By the finite-to-one nature of $\pi$ we may lift this 
ball-to-ball collision to any cell $P_{g_1}\times P_{g_2}$ of $N$, 
nevertheless, this choice determines the lift for the whole trajectory uniquely. 
Let us choose $g_1=g_2=Id$, as a consequence, the collision is lifted to $C_e$. 
By statement 1) of Corollary~\ref{welldef}, the same holds for all ball-to-ball collisions 
that 
belong to the island containing the fixed collision. 

Case (1) is a straightforward consequence of statement 2) from Corollary~\ref{welldef}.

In case (2) we may distinguish odd and even islands of the (semi)tra\-jectory 
depending on their ``distance'' from the ball-to-ball collision fixed above. 
More precisely, 
to define the parity of the island, consider the long collision sequence 
between the fixed ball-to-ball collision and any ball-to-ball collision 
of the island, and count the non-simple elements $\g_i$ on it. By 
straightforward application of Corollary~\ref{welldef}, statements 1) and 3)
any collision in an even island is lifted to $C_e$. On the other hand, by 
the same statements, there is a unique cylinder to which all collisions of 
odd islands are lifted. 
By Lemma~\ref{where} this unique cylinder is $C_R$ for some 
reflection $R$. 
\end{proof}
 
\begin{lemma}
\label{avoiding}
For the three polygons different from the square, the set of 
R-poor trajectories is slim. For all the four polygons the set of O-poor trajectories is slim 
and the set of R-poor, R$^+$-poor and O$^+$-poor trajectories has $\mu$ measure 0.
\end{lemma}

\begin{proof}
By Lemma~\ref{lemma1} it is enough to prove the statements in $N$.
We will use the upstairs characterizations of various types of poor points 
given by Lemma~\ref{poorup} and the description of the relevant cylinders
from Lemma~\ref{cylgeom}.

First consider the O-poor points in case $P$ is not 
the square. By Lemmas~\ref{poorup} 
and~\ref{cylgeom} the dynamics of such points is governed by a cylindric 
billiard which is (i) transitive (cf. Conjecture~\ref{trans}) and for which 
(ii) any two cylinders have transversal generator spaces. Theorem 2.4 from 
\cite{bcyl} states that such a system is a mixing semi-dispersing billiard. 
Furthermore O-poor points avoid all the
cylinders not of the form $C_g$, $g\in G_O$, an open set, thus applying the 
strong ball 
avoiding
theorem (Theorem~\ref{sballav}) yields the result.


In the case of the square the appropriate cylindric billiard is not 
transitive, there is an integral 
of motion, namely the projection of the velocity 
vector onto any of the two orthogonal base subspaces.  
We will prove that for each fixed value of this 
integral, O-poor points are slim on the surface of constant integral and 
thus 
by Lemma~\ref{intdim} are slim in the whole phase space. To see this, 
notice that on a
constant energy surface the system restricted to the set of 
O-poor points can be 
interpreted as
a direct product of two dispersing billiards which is mixing.  The open 
sets avoided have open intersections with
each constant energy surface, thus we can apply the strong ball avoiding 
theorem (Theorem \ref{sballav}) to conclude the case of O-poor points in the square.

In the R-poor case the appriopriate cylindric billiard consists of two 
``orthogonal'' cylinders. Dynamics, as discussed in \cite{S}, is a direct 
product of
a mixing semi-dispersing billiard system with a linear flow on $\S^1$. 
There is an integral of motion, 
namely the projection
of the velocity vector onto the direction of the linear flow,
i.e.~onto the common generator of the two cylinders.

This is the point where we should take into account that the geometry of the square is different 
from the other three polygons (cf. Remark~\ref{square}). 
Consider namely the direction of the above described linear flow factor and let us denote the
line of $\R^4$ parallel to it with $F$. In the non-square cases 
there is at least one of the avoided cylinders  
(e.g. any cylinder $C_O$ where $O$ is a rotation with a degree different from $180^\circ$) that has 
generator space not orthogonal to $F$. As a consequence, we may fix $p\in \S^1$ (in the direction 
of the linear flow) arbitrarily, the
erased cylinders intersect the leaf $\{p\} \times \T^{3}$ in some nonempty open set $U_p\subset \T^3$.
On the contrary, if the polygon $P$ is the square, all avoided cylinders 
have generator spaces orthogonal to $F$, and thus for certain points $p\in\S^1$ 
the intersection $U_p$ is empty. 

We may again prove slimness for each surface of constant integral and then integrate by 
Lemma~\ref{intdim}. Let us start with 
fixing the velocity component in the direction of the linear flow to be $0$, in such a situation we 
may disregard the 
integrable component. As to the mixing one, in the non-square cases the observations above show that 
there is an open set avoided by the whole trajectory. The strong ball avoiding theorem gives slimness 
of R-poor points within the surface. To see this observe that R-poor points 
necessarily collide with both cylinders of the orthogonal cylindric billiard 
upstairs (otherwise they would be O-poor) thus they belong to the full 
measure invariant set to which Theorem~\ref{sballav} applies (see also 
Remark~\ref{strongball}).
Note that the argument does not work for the 
square, however, points having zero velocity component in the direction of the 
linear flow form 
a set of zero measure.        

For the values of non-zero component in the integrable direction we do not need to treat the square 
separately. 
Fix $p \in \S^1$ for which $U_p$ above is nonempty 
(in the non-square cases $p$ can be chosen as any point of 
$\S^1$).  
Then there is an open neighborhood $U\subset \S^1$ of 
$p$, such that for
any $q \in U$,   $\{q\} \times \T^{3}$ intersects an erased cylinder in 
an open set.
Let $H$ be the set of times when the linear flow starting from the 
projection of the point $x$
visits $U$.  This set is syndetic, i.e.~it has bounded gaps.
We apply the strong ball avoiding theorem (Theorem \ref{sballav}) to the mixing component 
(i.e.~integral surface)
of the system, and integrate over the integrals of motion by Lemma~\ref{intdim} to conclude.

The reasoning above shows that the set of R-poor points is of zero measure, and that it is 
slim except for the case of the square.

Finally we turn to O$^+$-poor and R$^+$-poor points.  The arguments are identical
to the above ones
up to three formal changes:
trajectories are replaced by semitrajectories, the strong ball avoiding 
theorem
is replaced by the weak ball avoiding theorem, and the integration is done via 
Fubini's theorem.
\end{proof}

\begin{corollary}\label{noname}
The set of $x \in M$ such that $x$ is not hyperbolic has $\mu$ measure 
$0$.
\end{corollary}

\begin{proof}
It is enough to consider regular phase points ($x\in M^0$) 
since their complement 
has $\mu$-measure $0$, cf. Lemma~\ref{singslim}.


Recall Definitions~\ref{rich} (richness) and \ref{suff} (sufficiency).
In Lemma~\ref{hypgeom} we have shown that almost all rich phase points
are sufficient. The local hyperbolicity theorem (Theorem~\ref{lochyp}) 
states that every 
sufficient 
phase point has a neighborhood in which almost every point is hyperbolic.
Thus we only need to prove that almost every point is rich.

Unextendable, R-poor, O-poor, R$^+$-poor, and O$^+$-poor are all of measure
0.  If a phase point belongs to the complement of all of these then it is
rich.
\end{proof}

The following Lemma, on the one hand, exploits transversality in a way 
analogous to Lemma~\ref{transv1} and, on the other hand, is a strengthening of 
Lemma~\ref{avoiding} for the case of poor semitrajectories.

\begin{lemma}
\label{transv2} 
Consider the one-codimensional submanifold $L$ defined by the relation that 
$v_1(t_0)-v_2(t_0)$ is restricted to a line. The set of points that belong to 
$L$ and are
R$^+$-poor or O$^+$-poor is slim. 
\end{lemma}
\begin{proof}
Just like in the proof of Lemma~\ref{avoiding}, we work in the upstairs 
phase space $N$. The measure induced by $\nu$ on (the lift of) $L$ will be 
referred to as $\nu_L$.
Let us consider first the O$^+$-poor case for a polygon different from the square. 
Such points necessarily 
avoid an open set in positive time, the inner radius of which we will denote by 
$2\varepsilon$. 
We erase the cylinders constituting the avoided open set. 
By the geometry of the cylinders not erased the 
corresponding mixing semi--dispersing billiard dynamics is mixing. 
Unlike Lemma~\ref{avoiding}, in all arguments in the rest of the proof 
we do not consider only the set of 
O$^+$-poor points -- where the original and the modified dynamics coincide 
-- but the full modified billiard system.

As the set of \op points is closed (cf Definition~\ref{poor}), 
we need to prove slimness for a closed set and thus may use the 
characterization of
Lemma~\ref{adddim}. Assume the contrary: 
the set of O$^+$-poor points has a nonempty interior $A$ inside $L$. By the 
hyperbolicity of the modified billiard, there is 
$\hat{A}\subset A$ with
$\nu_L(\hat{A})=\nu_L(A)$ such that any $y\in \hat{A}$ has a
local stable 
manifold 
$\gamma^s_{\varepsilon}(y)$ of 
some positive length, fixed to be less than $\varepsilon$. As the modified billiard is mixing, 
these manifolds are 
strictly concave local orthogonal manifolds, thus they are transversal to 
$L$ (cf. the proof of Lemma~\ref{transv1}). Thus 
\be
\label{furry}
\hat{B}:= \bigcup_{y\in \hat{A}} \gamma^s_{\varepsilon}(y)
\ee  
is a set of positive $\mu$-measure. On the other hand the points of $\hat{B}$ 
avoid an open set (of inner radius $\varepsilon$), which is, by the weak 
ball-avoiding theorem (Theorem~\ref{ballav}) a contradiction.

As to the R$^+$-poor case (and the \op case for the square) 
the argument is analogous with the following modifications. 
First we take a finite cover of the set of R$^+$-poor points indexed by the 
reflection $R\in G$, where the semitrajectory collides only with $C_e$ and
$C_R$. 
We show that the intersection of $L$ with any element of this finite cover 
is slim. (This first step is irrelevant for the \op case in the square.)
The sets $A$ and 
$\hat{A}$ are constructed as above, however, this time the stable manifolds 
are not necessarily strictly concave in all directions. Nevertheless, they 
are strictly concave in the scattering direction of the first ball-to-ball 
collision after time moment $t_0$. Thus they are transversal to $L$.  
\end{proof}

\begin{lemma}
\label{regslim}
For any polygon different from the square, 
the set of phase points that are regular and non-sufficient is slim.
\end{lemma}

\begin{proof}
As the polygon $P$ is different from the square, 
the union of non-extendable, R-poor and O-poor points is a slim set 
(cf. Lemma~\ref{avoiding}), thus we 
restrict to its complement. The rest of the proof holds true for 
all the four polygons.  

In the characterization below rotations and 
reflections are arbitrary unless the axis of the reflection is specified.
One of the following cases applies:
\begin{enumerate}
\item the point is \op and there is a finite segment with short collision 
sequence $b\,\g_1\,(b,s)\,\g_2\,b$, where $\hg_1=R$, $\hg_2=O$ and $t_0$ is 
a time moment just before the 
island $(b,s)$.
\item the point is R$^+$-poor and there is a finite segment with short 
collision 
sequence $b\,\g_1\,(b,s)\,\g_2\,b$, where $\hg_1=R_E$, $\bg_2=R_{E'}$ with 
$E'\ne sE$, and $t_0$ is a time moment just before the 
island $(b,s)$.
\item the point is R$^+$-poor and there is a finite segment with short 
collision 
sequence $b\,\g_1\,(b,s_1)\,\g_2\,(b,s_2)\, \g_3 b$, where $\hg_1=R'$, 
$\hg_2=O$, $\hg_3=R$ ($R=R'$ not excluded) and $t_0$ is a time moment just 
before the 
island $(b,s_2)$.
\item the point is R$^+$-poor and there is a finite segment with short 
collision sequence 
$b\,\g_1\,(b,s_1)\,\g_2\, (b,s_2) \dots 
(b,s_{k-2}) \, \g_{k-1} (b,s_{k-1}) \g_k b$,\\ where $\hg_1=R'$, $\hg_2=O'$, 
$\hg_{k-1}=O$, $\hg_k=R$ (neither $R'=R$, nor $O'=O$ excluded) and  
$t_0$ is a time moment just before the island $(b,s_{k-1})$.
\item the point is \rp and O$^-$-poor at the ``same time''. More precisely, 
there exists a finite segment with short symbolic sequence $b\,\g_1\,(b,s)\,\g_2\,b$, 
and time moments 
$t_0$ and $t_0'$ just before and after the island $(b,s)$, 
respectively, such that the semitrajectories starting at $t_0$ 
and ending at $t_0'$ have the relevant poorness properties.
\item the point is neither \op nor \rp.
\end{enumerate}

Before turning to the particular cases we note that for any \rp or \op point 
and non-collision time $t_1>t_0$ just preceding a ball-to-ball collision
the semitrajectory starting at $t_1$ can be used to conclude that 
the point is \rp or \op.    

We start by the observation that case (5), when considered in backward time, 
reduces to case (1).

As to case (1) we introduce one more non-collision time moment $t_1$ just 
before the last ball-to-ball collision on the sequence. From Lemma~\ref{hypgeom} 
(with respect to part (2) of the definition of rich)
we know that the point is sufficient unless either $v_1(t_0)-v_2(t_0)$ or 
$v_1(t_1)-v_2(t_1)$ is restricted to a line. However, by Lemma~\ref{transv2}  
the set of those \op points for which any of these relations apply is 
slim. 

In case (2) the argument is identical with reference to Lemma~\ref{hypgeom} 
part (1) this time.

In case (3) we introduce, in addition to $t_0$, non-collision time moments 
$t_1$,$t_2$ and $t_3$ just before, just after the island $(b,s_1)$, and just after the 
island $(b,s_2)$, respectively. We apply Lemma~\ref{hypgeom}, 
(with respect to part (2) of the definition of rich) to the 
subsequence $b\,\g_1\,(b,s_1)\,\g_2\,b$: nonsufficiency is only possible if 
$v_1-v_2$, either 
at $t_0$ or at $t_1$ is restricted to a line. If restriction appears at $t_0$ 
slimness follows from Lemma~\ref{transv2}. Otherwise we apply 
Lemma~\ref{hypgeom}, part (2) to the sequence $b\,\g_2\,(b,s_2)\,\g_3\,b$ in 
backward time: 
sufficiency appears unless $v_1-v_2$ is restricted to a line either at $t_2$ 
or at $t_3$ (see Remark~\ref{firstcodim}). 
Any of these restrictions gives, together with the one at $t_1$, 
two transversal codimensions in the sense of Lemma~\ref{transv1}. To see this 
observe there is at least one ball-to-ball collision 
between $t_1$ and $t_2$. Thus we have 
sufficiency apart from a slim set in case (3).

Finally, points belonging to either of the cases (4) or (6) are twice rich 
thus, apart from a slim subset, are sufficient (cf. Corollary~\ref{2rich}).
\end{proof}

Finally we turn to the analysis of singular trajectories.  We extend
the definitions of extendable, \op and \rp to singular trajectories 
with the only additional restriction that we need $t_0 > 0$.

\begin{lemma}
\label{sing}
Within the set $S^+ \cap M^1$, the set of \rp and \op points
has 0 measure with respect to $m_{S^+}$ (the induced Riemannian measure 
on $S^+$, cf. Subsection~\ref{sec2.3}). 
\end{lemma}

\begin{proof}
The proof of this lemma is, on the one hand, rather standard (it is a 
straightforward adaptation of the analogous statements from the literature, 
e.g.\ Lemma 6.1 in \cite{inv}, part I) and, on the other hand, similar to the 
proof of Lemma~\ref{transv2} in this paper. Thus we only give a sketch.
We begin by considering the \op case for $P$ different from the square. 
Assume the contrary of the statement: 
there is
a subset $A$ of $S^+ \cap M^1$ of positive $m_{S^+}$-measure such that 
any $x \in A$ is \op. This implies that the semitrajectory 
$\{\Psi_t x: t \ge t_0\}$ (i) coincides with that of a modified 
billiard dynamics, (ii) avoids
an open set (the erased cylinders). The modified semi-dispersing billiard is 
mixing and hyperbolic, thus for almost all points $y\in A$ there exist local 
stable manifolds of positive inner radius at $\Psi^{t_0}y$ which we denote as 
$\gamma^s(\Psi^{t_0}y)$. These stable manifolds are (strictly) concave 
orthogonal manifolds, thus (by Lemma 4.8 in \cite{KSSz}) their pre-images are 
transversal to $S^+$. We may construct a set analogous to (\ref{furry}) that 
has, on the one hand, by the above mentioned transversality positive measure and avoids, on the other hand, an open set and thus has zero measure by the weak 
ball avoiding theorem (Theorem~\ref{ballav}). Thus we get a contradiction.  
 
The \rp case (and furthermore, the \op case for the square) is similar with 
the only difference that the modified dynamics is not fully hyperbolic: the 
local stable manifolds are 2 dimensional. To obtain three dimensional concave 
orthogonal manifolds we combine the stable manifolds with infinitesimal lines 
in the neutral direction of the modified billiard dynamics 
(on details see, e.g.  Lemma 6.1 in \cite{inv}, part I).
\end{proof}

This Lemma has two important immediate consequences.

\begin{corollary}
\label{ansatz}
(1) The Sinai-Chernov Ansatz holds for our system.

(2) The set of non-sufficient points that belong to $M^1$ is slim. 
\end{corollary}

\begin{proof}
To prove (2) we may assume $x\in M^1$ to be neither \op nor \rp. 
To see this note that the complement is slim by 
Lemma~\ref{sing}, as sets of zero volume measure 
must have empty interior and the characterization of Lemma~\ref{adddim} applies.
This however means that $x$ is twice rich and thus, apart from a slim subset, sufficient by 
Corollary~\ref{2rich}.

To prove (1) we restrict our attention to $x\in S^+\cap M^1$ as the rest is of 
zero measure in $S^+$ (cf. Lemma~\ref{singslim}). 
Furthermore, by Lemma~\ref{sing} we may assume $x$ to 
be neither \op nor \rp. Thus there certainly exists a subsegment of the 
semitrajectory $\{\Psi_t x: t \ge t_0\}$ with collision sequence 
$b\,\g_1\,(b,s)\,\g_2\,b$ where $\bg_2=R_E$ and either $\hg_1=O$ or 
$\hg_1=R_{E'}$ with $E'\ne s^{-1}E$. 
We apply Lemma~\ref{hypgeom} and Remark~\ref{firstcodim} in {\it 
backward time}: $x$ is 
sufficient unless it belongs to a one-codimensional submanifold $L$ defined by 
the condition that $v_1(t)-v_2(t)$ is restricted to a line 
for some $t(>t_0>0)$ {\it just after} some ball-to-ball collision. To prove 
that non-sufficient points are of zero measure in $S^+$ it is enough to show 
the transversality of $L$ to $S^+$. To see this we apply the strategy of 
Lemma~\ref{transv1}, we foliate $L$ with the equivalence classes of the 
relation (\ref{eqclass}). The pre-images of these equivalence classes are 
(strictly) concave local 
orthogonal manifolds and thus transversal to $S^+$.  
\end{proof}

Now we can easily prove the main theorem of the paper.

\begin{theorem}
\label{main} If $R$ is sufficiently small so that the phase space is connected
then 
the system of two balls in any of the four integrable polygons is ergodic.
Otherwise the system is ergodic on each connected component of the phase space.
\end{theorem}

\begin{proof}
As the boundary of the billiard table is defined by algebraic equations the 
semi-dispersing billiard is algebraic. In addition, by Corollary~\ref{ansatz} the 
Sinai-Chernov Ansatz holds. Thus the Local Ergodicity Theorem applies. 

To conclude that the local ergodic components make up a single ergodic component 
it is enough to show that the set of non-sufficient points is slim. This is, in fact, 
not true for the square, for the case of which, however, ergodicity 
(in contrast to hyperbolicity proved in Corollary~\ref{noname}) easily follows from 
\cite{S} (cf. Remark~\ref{square}).

Consider the other three polygons. We may 
restrict to $M^1\cup M^0$ as the complement is slim. However, by 
Lemma~\ref{regslim} and Corollary~\ref{ansatz}, respectively,
non-sufficient points of both $M^0$ and $M^1$ belong to slim subsets of $M$.
\end{proof}

\begin{remark}\label{noname2}
By \cite{KS,CH,OW} our systems are automatically K-mixing and Bernoulli.
\end{remark}

\section{Further results and outlook}
\setcounter{equation}{0}

In this final section we describe three problems in decreasing order of 
interest and difficulty. 

\subsection{Higher genus cases}\label{noname5} 
A natural question that arises is if one can 
extend our result to other polygons, in particular to rational ones, i.e.\ 
those for which all angles between sides are rational multiples of $180^\circ$. 
The one point particle dynamics in rational polygons essentially reduces to 
the study of the linear flow on a flat surface with singularities. The case
of genus 1 (the flat torus without any singularity) corresponds to the four
polygons treated in this article. It would be very interesting to extend  
our result in the higher genus case. There are various additional 
difficulties, mainly because in place of cylindric billiards on 
$\T^4(=\T^2 \times \T^2)$ one should analyze the flows on manifolds 
$S \times S$ with cylindric 
regions removed, where $S$ is a higher genus flat surface. 

\subsection{Higher dimensional case}\label{noname3} With the techniques of our 
article one can prove that the motion of two hard balls in
any right prism of the form $P \times [0,1]^n$, where $P$ is one of the 
four polygons of Theorem~\ref{main}, is hyperbolic and ergodic. Just like in 
the article, the system can be lifted up to a cylindric billiard on 
$\T^{2n+2}$.
However, more types of 
rich trajectories arise (except in the cube) since the group $G$ is more 
complicated. Namely, analogously to reflections and rotations, one most 
distinguish group elements with or without non-trivial invariant subsets. 
These later ones, unlike the two dimensional case of rotations, do not form  
a subgroup of $G$. This problem does not appear in the cube, which has 
remarkable 
symmetries: the group $G$ is commutative, the cylindric billiard upstairs is 
orthogonal (in the sense of \cite{S}) and even though the number 
of rich sequences increase because of higher dimensionality, they all belong, 
essentially, to one of the types discussed in our paper. We choose not to 
analyze the
issue of right prisms here 
to keep the length of the article reasonable.

\subsection{Invisible corners}\label{noname4}
\begin{figure}[hbt]
\centering\resizebox{12cm}{!}{
\subfigure{\epsfig{file=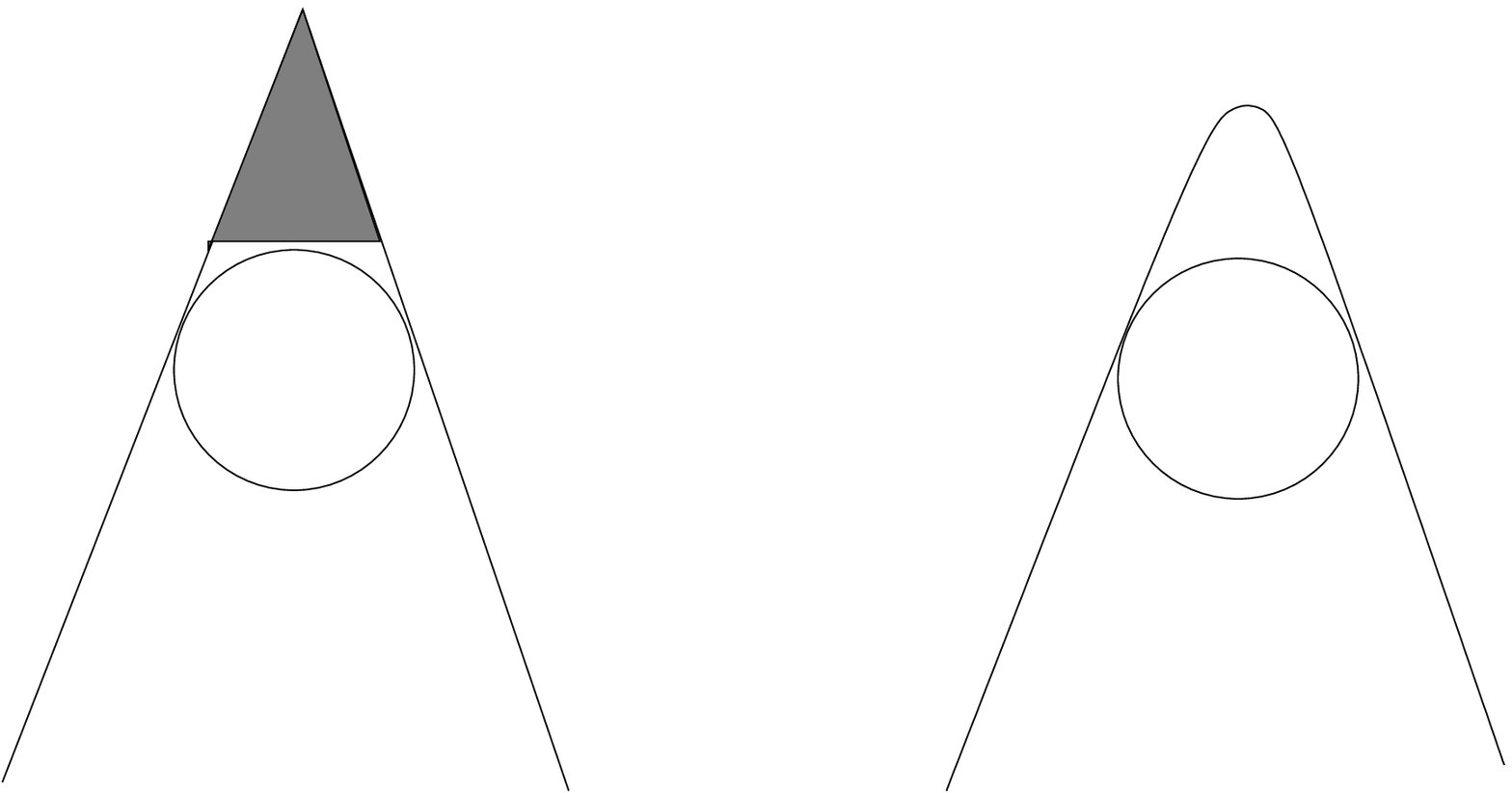}}\qquad\qquad\qquad
\subfigure{\epsfig{file=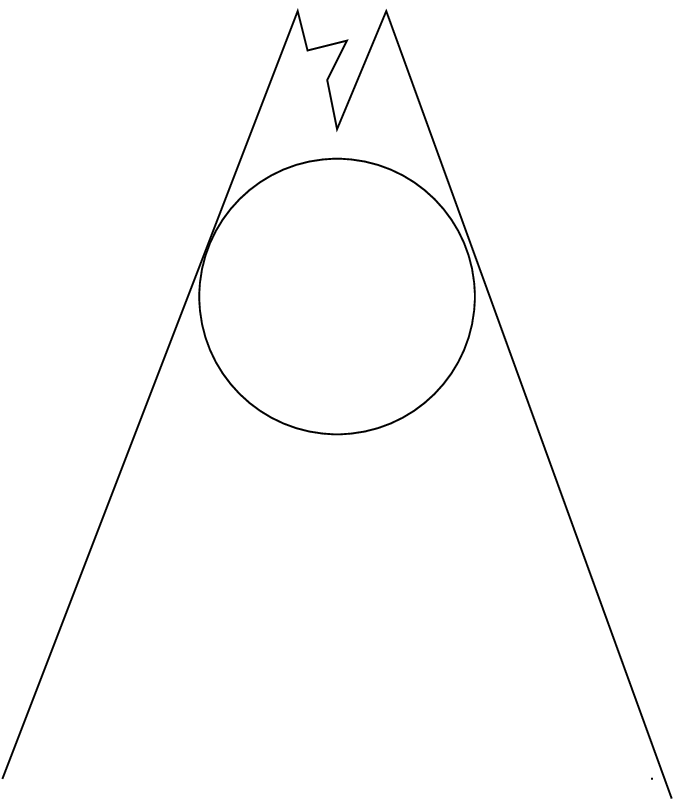}} } 
\caption{(a) Invisible corners\quad (b) smooth table and\quad (c) fractal table} 
\label{corner}
\end{figure}

Fix a convex polygon $P$ and the (common) radius of the 
balls $r$. Place one of the balls in a corner as in Figure~\ref{corner} and 
consider the shaded region blocked by the ball. Note that neither of the balls 
ever collides with the part of the boundary of $P$ that intersects the shaded 
region. Thus we can modify this part of the boundary in any way that 
stays within the shaded region and 
this will not effect the dynamics at all. In particular if either the one 
or the two ball
system in $P$ is ergodic, it remains ergodic for the modified table. The 
modification may result in a rational or irrational polygon, an infinite 
polygon, a Sinai or a Bunimovich billiard table or even a fractal. This way 
strange results may arise: e.g. both for the one and the two ball cases 
we can construct ergodic $C^\infty$ tables, see Figure~\ref{corner}.

\section*{Acknowledgements}

P.B. is grateful for the hospitality and the inspiring research atmosphere of 
CPT Marseille. The financial support of the Hungarian National Foundation 
for Scientific Research (OTKA), grants T32022 and TS040719 is thankfully acknowledged.

\end{document}